# Approximation of Metric Spaces by Reeb Graphs: Cycle Rank of a Reeb Graph, the Co-rank of the Fundamental Group, and Large Components of Level Sets on Riemannian Manifolds

Irina Gelbukh

*CIC, Instituto Politécnico Nacional, 07738, Mexico City, Mexico*

**Abstract.** For a connected locally path-connected topological space $X$ and a continuous function $f$ on it such that its Reeb graph $R_f$ is a finite topological graph, we show that the cycle rank of $R_f$, i.e., the first Betti number $b_1(R_f)$, in computational geometry called *number of loops*, is bounded from above by the co-rank of the fundamental group $\pi_1(X)$, the condition of local path-connectedness being important since generally $b_1(R_f)$ can even exceed $b_1(X)$. We give some practical methods for calculating the co-rank of $\pi_1(X)$ and a closely related value, the isotropy index. We apply our bound to improve upper bounds on the distortion of the Reeb quotient map, and thus on the Gromov-Hausdorff approximation of the space by Reeb graphs, for the distance function on a compact geodesic space and for a simple Morse function on a closed Riemannian manifold. This distortion is bounded from below by what we call the *Reeb width* $b(M)$ of a metric space $M$, which guarantees that any real-valued continuous function on $M$ has large enough contour (connected component of a level set). We show that for a Riemannian manifold, $b(M)$ is non-zero and give a lower bound on it in terms of characteristics of the manifold. In particular, we show that any real-valued continuous function on a closed Euclidean unit ball $E$ of dimension at least two has a contour $C$ with $\mathrm{diam}(C \cap \partial E) \geq \sqrt{3}$.

## 1. Introduction

Given a topological space $X$, the Reeb graph $R_f$ of a continuous function $f : X \to \mathbb{R}$ is a space obtained by contracting the connected components of the level sets of $f$ to points, endowed with the quotient topology (see Section 2.1). It shows the evolution of the topology of the level sets, thus providing important information on the behavior of the function. The Reeb graph was introduced in 1946 for Morse functions on a compact manifold in the context of Morse theory [39]. Later it was used for study and classification of more general types of functions, e.g., Morse-Bott functions or functions with isolated singularities on a compact manifold [22, 27, 40].

The cycle rank, defined as the first Betti number $b_1(R_f)$, is an important invariant of the Reeb graph. It reflects the complexity of the graph and appears in various bounds involving Reeb graphs [8, 30, 31]. Michalak [33] defined the *Reeb number* $\mathcal{R}(M)$ of a closed manifold $M$ (see end of Section 3) as the maximum cycle rank (which he calls *number of cycles*) of the Reeb graph of a smooth function with finitely many critical points on this manifold. He uses $\mathcal{R}(M)$ to characterize the set of graphs that can be realized as the Reeb graphs of such functions.

In the 1990s, the notion of Reeb graph was generalized, under the name of *foliation graph*, from Morse functions to Morse forms (closed 1-forms with Morse singularities) [28] and was extensively used to study the foliations defined by such forms. Its cycle rank coincides with the number of homologically independent compact leaves of the foliation [16]. This number gives a lower bound on the number of conic singularities of the form [17], as well as an upper bound on the number of minimal components, which, in particular, gives a condition for the foliation to have no minimal components [15]. In the latter case, the cycle rank of the foliation graph gives an upper bound on the rank of the form [16].

The notion of Reeb graph also has important applications in computational geometry [4], where its cycle rank $b_1(R_f)$ is known under the name of the *number of loops of the Reeb graph* [9, 20]. In machine learning and big data analysis, graphs equipped with a metric derived from the data represent hidden structure in complex data sets [14]; Reeb-type graphs are useful for approximating metric spaces, with the cycle rank being used to bound such approximation [8, 30, 31] (see Section 5).

As a general quotient space, the Reeb graph can be quite ill-behaved: e.g., it may be non-Hausdorff or even not one-dimensional topological space. The majority of existing studies are restricted to functions on compact manifolds with some conditions of finiteness on the critical set, such as Morse or Morse-Bott functions or functions with isolated critical points, under which the Reeb graph is a finite topological graph. Polulyakh [36, 37] extends the study of Reeb graphs to non-compact surfaces, studying conditions on a continuous function on such a surface under which its Reeb graph has a particularly simple structure, being a locally finite one-dimensional space called *graph with stalks*. In this paper, we further consider general topological spaces, not necessarily manifolds. However, we restrict our study to those continuous functions whose Reeb graphs are finite topological graphs.

The cycle rank of the Reeb graph of a continuous function is widely believed to be bounded from above by the first Betti number of the domain $X$ of the function:

$$b_1(R_f) \leq b_1(X). \tag{1}$$

We show, however, that this "obvious" bound does not generally hold; this misconception has led some authors to generally wrong statements (see Section 3). In fact, for even not too ill-behaved spaces and quite well-behaved functions, the values $b_1(R_f)$ and $b_1(X)$ are completely unrelated (Proposition 3.2).

On the other hand, for well-behaved spaces and functions, including those typically used in practical applications, we show a much stronger estimate. Namely, for a connected locally path-connected topological space $X$ and continuous function $f : X \to \mathbb{R}$ whose Reeb graph $R_f$ is a finite topological graph, we show that the the cycle rank $b_1(R_f)$ is bounded by the co-rank of the fundamental group of $X$ (Theorem 3.1):

$$b_1(R_f) \leq b_1'(X) \stackrel{\text{def}}{=} \operatorname{corank}(\pi_1(X)). \tag{2}$$

This bound is tight in the class of smooth closed manifolds (Proposition 3.9).

We give some practical methods for calculation of the value $b_1'(X)$ for spaces composed from simpler ones (Theorem 4.1). In turn, this value can be bound in terms of the so-called isotropy index $h(X)$—the maximum rank of a subgroup in $H^1(X; \mathbb{Z})$ with trivial cup product; we also calculate $h(X)$ for spaces composed from simpler ones (Theorem 4.3) Some of these results have been previously reported for smooth (orientable) manifolds; here we extend them to more general topological spaces.

As an application of our bound, we improve three estimates on the approximation of some compact metric spaces by graphs. A compact connected length space $X$ can be approximated by finite metric graphs $G$ under the Gromov–Hausdorff distance $d_{GH}(X, G)$ [6, Proposition 7.5.5] (see Section 2.4). While an arbitrary such graph can be quite complex, the well-studied Reeb graphs of continuous functions are useful for this purpose due to their simple structure. Mémoli and Okutan [31] studied how well a compact geodesic space $(X, d)$ can be approximated by graphs $G$ with low enough $b_1(G)$ by bounding $d_{GH}(X, G)$ from below and from above in terms of $d_{GH}(X, R_f)$ and $b_1(R_f)$, where $R_f$ is the Reeb graph for the distance function $f = d(p, \cdot)$, $p \in X$, under the assumption that $R_f$ is a finite topological graph. To estimate $b_1(R_f)$, the authors used the trivial bound (1). We show that using our bound (2) instead improves their result (Section 5.1), in some cases by a large margin (Example 5.1).

Moreover, various works have been devoted to further bound the distance $d_{GH}(X, R_f)$ from above for selected types of spaces and functions [8, 30, 43]. It is known [6, Theorem 7.3.25] that $2d_{GH}(X, R_f) \leq \text{dis}(\varphi)$, the distortion of the Reeb quotient map $\varphi: X \to R_f$ (Sections 2.1, 2.4). Chazal et al. [8] bounded $\text{dis}(\varphi)$ in terms of $R_f$ for the function $f$ from [31] described above, while Mémoli and Okutan [30] bounded $\text{dis}(\varphi)$ in terms of $b_1(R_f)$ for a simple Morse function $f$ on a closed Riemannian manifold $M$. To further estimate the distortion $\text{dis}(\varphi)$, and thus $d_{GH}(X, R_f)$, in terms of characteristics of the space $X$, the authors of both papers used the trivial bound (1). Our bound (2) improves their results (Sections 5.2, 5.3), giving much stronger estimates (Theorems 5.2, 5.4) for some important spaces (Examples 5.3, 5.5, 5.6).

While the upper bounds guarantee the possibility of approximating the space by Reeb graphs with certain quality, a lower bound independent from the function marks the limits on the best possible approximation. We show (Proposition 5.8) that for a space containing a homeomorphic image of $\mathbb{R}^2$, the lower bound on $\text{dis}(\varphi)$ is positive; namely, we show that every real-valued continuous function on such a space will have large enough contours.

We introduce the *Reeb width* of a metric space $X$ (Definition 6.1), which guarantees that for any real-valued continuous function on $X$, some contours will be sufficiently large:

$$b(X) = \inf_{\substack{\text{continuous} \\ f: X \to \mathbb{R}}} \sup_{\text{contours } C \text{ of } f} \text{diam}(C);$$

obviously, for the distortion of the Reeb map of any continuous function on $X$, it holds $\text{dis}(\varphi) \geq b(X)$. The notion of Reeb width is closely related with that of one-dimensional Urysohn (Alexandrov) width $u_1(X)$, but studying this relationship is left for future work.

For Riemannian manifolds, we give specific lower bounds on their Reeb width (roughly speaking, the diameter of the largest contour of any real-valued continuous function) in terms of local (Proposition 6.4) and global (Theorem 6.7, Corollary 6.8) characteristics of the manifold. The bound in terms of local characteristics is important because in some cases the global characteristics provide too coarse information (Examples 6.9, 6.10). For some spaces this bound is sharp (Example 6.6).

Finally, we refine an important result of Maliszewski and Szyszkowski [26] about diameters of contours on disks. They have shown that any real-valued continuous function $f$ on a unit disk $E = E^n$ has a contour $C$ with $\text{diam}(C) = 2$ for $n \geq 3$ and $\text{diam}(C) \geq \sqrt{3}$ for $n = 2$. We show that $f$ will even have a contour $C$ with $\text{diam}(C \cap \partial E) = 2$ or $\text{diam}(C \cap \partial E) \geq \sqrt{3}$, respectively. This is trivial for $n \geq 3$, so we give a proof for $n = 2$. We also give a similar statement for an open disk.

The paper is organized as follows. In Section 2, we introduce some necessary notions and facts. In Section 3, we give our main results on the cycle rank of the finite Reeb graph. In Section 4, we give some practical methods for calculation of the co-rank of the fundamental group for topological spaces composed of simpler ones. In Section 5, we apply our bound on the cycle rank of the Reeb graph to improve three known bounds on the distortion of the Reeb quotient map. We also show that there is a non-zero lower bound on the distortion, which in Section 6 we calculate for Riemannian manifolds in terms of their local and global characteristics. We also show that every real-valued continuous function on a closed Euclidean ball has a contour whose intersection with the boundary of the ball has large diameter.

## 2. Definitions and Useful Facts

In this section, we give some definitions and facts useful for the description of our main results.

### 2.1. Reeb Graph and Simple Morse Functions

Given a topological space $X$, by a *contour* of a map $f: X \to Y$ we understand a connected component of the level set $f^{-1}(y)$ for some $y \in Y$.

The *Reeb graph* (called by some authors *Kronrod-Reeb graph* or *Reeb quotient space*) $R_f$ of a continuous function $f : X \to \mathbb{R}$ is the quotient space $X/\sim$, the equivalence relation $x \sim y$ holding whenever $x$ and $y$ belong to the same contour of $f$, endowed with the quotient topology. The quotient map $\varphi : X \to R_f$ is called the *Reeb quotient map*. The Reeb quotient map is continuous.

A smooth function on a smooth manifold is called a *Morse* function if all its critical points are non-degenerate. The Reeb graph defined by a Morse function on a compact manifold is a finite topological graph [39, Theorem 1].

A Morse function is called *simple* if each its critical level contains only one critical point; such functions are also called excellent [30] or non-resonant [34]. The Reeb graph defined by a simple Morse function has especially simple structure: its vertices have degree one or three. The set of simple Morse functions is dense in the space of continuous functions and open in the space of smooth functions [3, Theorem 5.31], so it is natural to use them for approximating.

**Lemma 2.1.** *Let $M$ be a connected smooth closed manifold $M$, and $N \subset M$ a two-sided compact codimension-one submanifold. Then connected components of $N$ are contours of some simple Morse function $f : M \to \mathbb{R}$.*

*Proof.* Consider a metric on $M$. Consider a product neighborhood $U = N \times [0, 1]$ with $N \times 0 = N$. Cutting $M$ open by $\partial U$ results in two manifolds with boundary, $\tilde{U}$ and $\tilde{M}$, such that $M$ is obtained by appropriately gluing $\partial \tilde{U}$ with $\partial \tilde{M}$.

On $\tilde{M}$, one can choose a simple Morse function $f_{\tilde{M}}$ close to the distance function from the boundary $\partial \tilde{M}$, equal to zero on $\partial \tilde{M}$, and having no critical points on $\partial \tilde{M}$. Similarly, consider such a function $f_{\tilde{U}}$ on $\tilde{U}$. Both functions can be chosen in such a way that the functions $f_{\tilde{M}}$ and $-f_{\tilde{U}}$ fit together smoothly into a simple Morse function $f$ on $M$, which has the desired properties. □

### 2.2. Co-rank of the Fundamental Group

The *co-rank* of a finitely generated group $G$ [23], also known as the inner rank [24], is the maximum rank of a free homomorphic image of $G$. For a path-connected topological space $X$, we denote $b'_1(X) = \text{corank}(\pi_1(X))$. This notation was introduced by Arnoux and Levitt [2] under the term "non-commutative Betti number;" in case of 3-manifolds this value is called the cut-number [41]. If the fundamental group $\pi_1(X)$ is finitely generated, as in the case of compact manifolds, then $b'_1(X)$ is finite. Obviously, $b'_1(X) \leq b_1(X)$. Some methods of calculating $b'_1(X)$ can be found in Section 4 below.

**Example 2.2.** *Denote by $T^n$ an n-torus, by $M^2_g = \#^g T^2$ the closed orientable surface of genus $g$, and by $N^2_g = \#^g \mathbb{R}P^2$ a closed non-orientable surface. Then*

$$b'_1(T^n) = 1 \quad [18], \qquad b_1(T^n) = n;$$
$$b'_1(M^2_g) = g \quad [23], \qquad b_1(M^2_g) = 2g;$$
$$b'_1(N^2_g) = \left[\tfrac{g}{2}\right] \quad [18, \text{Eq. (4.1)}], \qquad b_1(N^2_g) = g - 1.$$

*Since a punctured surface is homotopy equivalent to a wedge of circles, its fundamental group is free. Denoting by $M^2_{g,h}$ and $N^2_{g,h}$ the corresponding connected compact surfaces or with $h \geq 1$ boundary components, we have:*

$$b'_1(M^2_{g,h}) = b_1(M^2_{g,h}) = 2g + h - 1 \quad [9, \text{p. 346}];$$
$$b'_1(N^2_{g,h}) = b_1(N^2_{g,h}) = g + h - 1 \quad [9, \text{p. 347}].$$

*and similarly for the surfaces with $h$ points removed.*

**Proposition 2.3** ([10, Theorem 1], quoted by [32, Remark 5.3]). *Let $M$ be a connected closed smooth manifold. Then $b'_1(M)$ is the maximum number $k$ of disjoint compact connected codimension-one submanifolds $N_1, \ldots, N_k \subset M$ such that each $N_i$ has a product neighborhood and $M \setminus \bigcup_{i=1}^{k} N_i$ is connected.*

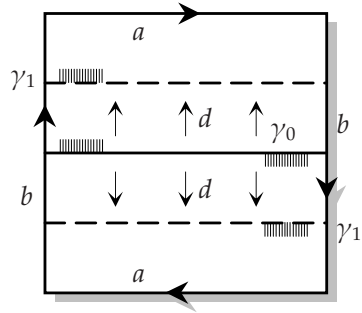

Figure 1: $P = \mathbb{R}P^2$, shown as a square with identified opposite sides *a* and *b* with the corresponding orientation. The hatching on the submanifold $\gamma_0$ (solid line) shows that there is no continuous normal field on it; the hatching at one side of the submanifold $\gamma_1$ (dotted line) shows that it has a product neighborhood.

**Remark 2.4.** *For a connected smooth closed manifold $M$ and a smooth function $f : M \to \mathbb{R}$, if $f$ has $k$ contours $N_1, \ldots, N_k$ such that $N_i \cap \mathrm{Crit}(f) = \emptyset$ for all $i$ and $M \setminus \bigcup_{i=1}^{k} N_i$ is connected, then $k \le b'_1(M)$. Indeed, the normal bundle defined by the gradient of $f$ is trivial; thus the tubular neighborhood of a regular contour is a product neighborhood.*

The requirement for the submanifolds to have a product neighborhood is essential:

**Example 2.5.** *For the real projective plane $P = \mathbb{R}P^2$, it holds $H_1(P) = \mathbb{Z}/2\mathbb{Z}$ and thus $b'_1(P) = b_1(P) = 0$. However, its middle line $\gamma_0$ is a submanifold that leaves $P \setminus \gamma_0$ connected; see Figure 1. Note that $\gamma_0$ does not have a product neighborhood, since its normal bundle is non-trivial. It is a contour of the distance function, $d(\gamma_0, \cdot) = 0$, but 0 is a non-regular value of this function. A contour $\gamma_1$ given by $d(\gamma_0, \cdot) = 1$ is regular, but $P \setminus \gamma_1$ is not connected. By Remark 2.4, no regular contour $\gamma$ of any continuous function on $P$ can leave $P \setminus \gamma$ connected.*

### 2.3. Isotropy Index

The *isotropy index* $h(X)$ of a topological space $X$ is the maximum rank of a subgroup in $H^1(X; \mathbb{Z})$ with trivial cup-product.

**Proposition 2.6.** *Let $X$ be a path-connected topological space with a finitely generated fundamental group $\pi_1(X)$. Then*

$$b'_1(X) \le h(X) \le b_1(X), \tag{3}$$

*where $b'_1(X) = \mathrm{corank}(\pi_1(X))$ and $b_1(X)$ is the first Betti number.*

While the proof can be obtained by generalizing [19, Proposition 39] to topological spaces, we give here a shorter self-contained proof.

*Proof.* Consider an epimorphism $\pi_1(X) \twoheadrightarrow F_b$ of the fundamental group onto a free group of $b = b'_1(X)$ generators and the corresponding map $\varphi : X \twoheadrightarrow V$, where $V = \bigvee_{i=1}^{b} S^1$. The induced map $\mathrm{Hom}(\pi_1(V), \mathbb{Z}) \hookrightarrow \mathrm{Hom}(\pi_1(X), \mathbb{Z})$ is injective. Since for any group $G$, homomorphisms $G \to \mathbb{Z}$ factor through its abelianization $G/[G, G]$, the latter rewrites as $\mathrm{Hom}(H_1(V), \mathbb{Z}) \hookrightarrow \mathrm{Hom}(H_1(X), \mathbb{Z})$. For connected $X$ and $V$, the universal coefficient theorem gives $\varphi^* : H^1(V; \mathbb{Z}) \hookrightarrow H^1(X; \mathbb{Z})$. Since $H^2(V; \mathbb{Z}) = 0$, the cup product on $\varphi^* H^1(V; \mathbb{Z}) \subseteq H^1(X; \mathbb{Z})$ is zero. We obtain $b = \mathrm{rank}(\varphi^* H^1(V; \mathbb{Z})) \le h(X)$. □

The isotropy index can also be easily calculated or bounded using vector space-based techniques. For a module $R$, denote by $h(X; R)$ the maximum rank of a submodule in $H^1(X; R)$ with trivial cup-product $\smile : H^1(X; R) \times H^1(X; R) \to H^2(X; R)$; note that $h(X) = h(X; \mathbb{Z})$. By [19, Lemma 7], for a space $X$ we can consider cohomology with coefficients in a field:

$$h(X) = h(X; \mathbb{Q}). \tag{4}$$

**Proposition 2.7.** *Under the conditions of Proposition 2.6, denote $b_1 = b_1(X)$ and $b_2 = b_2(X)$ the Betti numbers and $k = \dim \ker \smile$, where $\smile : H^1(X) \times H^1(X) \to H^2(X)$ is the cup product. Then*

$$\frac{b_1 + kb_2}{b_2 + 1} \leq h(X) \leq \frac{b_1 b_2 + k}{b_2 + 1};$$

*in particular, if $b_2(X) = 1$, then $h(X) = \frac{1}{2}(b_1 + k)$. If $\smile$ is surjective, then $h(X) \leq k + \frac{1}{2} + \left(\left(b_1 - k - \frac{1}{2}\right)^2 - 2b_2\right)^{\frac{1}{2}}$.*

Indeed, by (4), it is enough to apply [29, Propositions 1–3] to the vector spaces $H^i(X; \mathbb{Q})$. A similar fact for manifolds has been given as [19, Proposition 16]. That paper gives more details on calculation of the isotropy index for manifolds and on its geometric meaning.

### 2.4. Distortion of a Map and the Gromov–Hausdorff Distance

The *distortion* of a map $\varphi : X \to Y$ between metric spaces $(X, d_X)$ and $(Y, d_Y)$ is defined as

$$\mathrm{dis}(\varphi) = \sup_{x, x' \in X} |d_X(x, x') - d_Y(\varphi(x), \varphi(x'))|.$$

The *Hausdorff distance* between two subsets $A$ and $B$ of a metric space is

$$d_H(A, B) = \inf\{r > 0 \mid A \subseteq U_r(B) \text{ and } B \subseteq U_r(A)\},$$

where $U_r$ is the $r$-neighborhood of a set in the metric space [6, Definition 7.3.1].

The *Gromov–Hausdorff distance* $d_{GH}(X, Y)$ between metric spaces $X$ and $Y$ is the infimum of $r > 0$ for which there exist a metric space $Z$ and its subspaces $X'$ and $Y'$ isometric to $X$ and $Y$, respectively, such that $d_H(X', Y') < r$ [6, Definition 7.3.10]. The Gromov–Hausdorff distance measures how far two compact metric spaces are from being isometric.

If $\varphi$ is surjective, then for the Gromov-Hausdorff distance between $X$ and $Y$ it holds $d_{GH}(X, Y) \leq \frac{1}{2} \mathrm{dis}(\varphi)$ [6, Theorem 7.3.25]. In particular, for the Reeb quotient map $\varphi$, we have

$$d_{GH}(X, R_f) \leq \frac{1}{2} \mathrm{dis}(\varphi).$$

### 2.5. Thickness of a Function

The thickness of a function, introduced by Mémoli and Okutan [30, Section 5], indicates how the volume of the level sets is distributed with respect to their diameters. Namely, given a smooth function $f : M \to \mathbb{R}$ on an $n$-dimensional manifold $M$, the *thickness* $T_f$ of $f$ is the infimum of the $(n-1)$-thickness of its contours. Here, the *k-thickness* of a path-connected metric space $N$ is

$$T_N^k = \frac{\mu^k(N)}{(\mathrm{diam}(N))^k},$$

where $\mu^k(N)$ denotes the $k$-dimensional Hausdorff measure on $N$. For any smooth function $f$ on a surface, it holds $T_f \geq 1$ [30, Remark 5.4], whereas this does not hold for $n \geq 3$.

### 2.6. Level Sets and Riemannian Manifolds

**Proposition 2.8** ([26, Theorems 3.6, 3.8]). *Let $E^n \subset \mathbb{R}^n$, $n \geq 2$, be a closed unit ball and $f : E^n \to \mathbb{R}$ a continuous function. Then there exists a contour $C$ of $f$ with*

  (i) *$n = 2$: $\mathrm{diam}(C) \geq \sqrt{3}$,*
  (ii) *$n \geq 3$: $\mathrm{diam}(C) = 2$.*

The *k-dimensional Urysohn width* (called also Alexandrov width or Urysohn diameter $\text{diam}_k(X)$, and sometimes referred to as Uryson width [21]) of a metric space $X$, which, among equivalent definitions, is defined [26] as

$$u_k(X) = \inf_{\substack{Y^k \\ \varphi: X \to Y^k}} \sup_{\text{level sets } L \text{ of } \varphi} \text{diam}(L), \tag{5}$$

where $Y^k$ runs over all $k$-dimensional simplicial complexes and $\varphi$ over all continuous functions. We are interested in lower bounds on $u_1(X)$. While the Urysohn width of a Riemannian manifold $M$ has been mostly studied in the context of codimension one and of bounding it from above [21], some lower bounds are known. Obviously,

$$u_k(M) \geq u_{k+1}(M), \tag{6}$$

so lower bounds on higher-dimensional widths are useful for our case. G. Perelman has shown that for closed Riemannian manifolds with non-negative sectional curvature, the Urysohn width is related with the volume:

$$c(n) \text{vol}(M) \leq \prod_{k=0}^{n-1} u_k(M) \leq C(n) \text{vol}(M)$$

for some positive constants $c(n)$ and $C(n)$ (quoted by [5, Section 8.4]), which together with (6) implies

$$u_1(M) \geq \sqrt[n-1]{c(n) \frac{\text{vol}(M)}{\text{diam}(M)}}. \tag{7}$$

In the remainder of this subsection, $M$ denotes a Riemannian manifold.

By $\sec(p)$ we denote the *sectional curvature* (Gaussian curvature for $\dim M = 2$) at a point $p \in M$. For a subset $U \in M$ we define $\sec(U) = \sup_{p \in U} \sec(p)$.

By $S_K$ we denote the model surface of constant Gaussian curvature $K$: for $K = 0$ this is the Euclidean plane $\mathbb{R}^2$, for $K > 0$ a sphere $S^2$, for $K < 0$ a hyperbolic plane $H^2$.

The following fact is adapted from a part of the Cartan–Alexandrov–Toponogov comparison theorem for small triangles:

**Theorem 2.9 ([6, Theorem 6.5.6]).** *Let $M$ be a 2-dimensional Riemannian manifold. Let $U \subset M$ be such that for any $x, u \in U$ there exists a shortest path $[xy]$ in $M$ that lies entirely in $U$.[1] Let the Gaussian curvature $\sec(U) \leq K$. Let $a, b, c \in U$ and $a', b', c' \in S_k$, with the distances $d(a, c) = d(a', c')$, $d(b, c) = d(b', c')$ and the angles between the geodesic segments $\angle acb = \angle a'c'b'$ at $c'$. Then $d(a, b) \geq d(a', b')$.*

The *exponential map* $\exp_p : E \to M$, $E \subseteq T_x M$, where $M$ is a Riemannian manifold, is defined as $\exp_p(x) = g(1)$, where $g : \mathbb{R} \to M$ is a geodesic with $g(0) = p$, $g'(0) = x$.

The *injectivity radius* $\text{inj}(p)$ at $p \in M$ is the radius of the largest geodesic ball around $p$ on which $\exp_p$ is a diffeomorphism. We denote $\text{inj}(M) = \inf_{p \in M} \{\text{inj}(p)\}$.

A ball $B \subseteq M$ around $p \in M$ is *convex* if the distance $d(\cdot, p)$ is convex on $B$, i.e., has non-negative definite Hessian.

The *convexity radius* $\text{conv.rad}(p)$ at $p \in M$ is the largest radius of a ball $B$ around $p$ that is convex and any two points in $B$ are joined by a shortest geodesic segment lying in $B$ [35]. Note that $B$ satisfies the conditions on $U$ in Theorem 2.9. We denote $\text{conv.rad}(M) = \inf_{p \in M} \text{conv.rad}(p)$.

**Proposition 2.10 ([35, p. 259]).** *Let $M$ be a compact Riemannian manifold and $K = \sec(M)$. Then*

  (i) $K > 0$: $\text{conv.rad}(M) \geq \min\left\{\frac{1}{2} \text{inj}(M), \frac{\pi}{2\sqrt{K}}\right\}$,
  (ii) $K \leq 0$: $\text{conv.rad}(M) = \frac{1}{2} \text{inj}(M)$.

---

[1] In [6], the set $U$ is called *convex*, but their use of this term differs from our definition below, so we used [6, Definitions 1.1.1 ff., 2.1.10 ff., 3.6.5] for its interpretation.

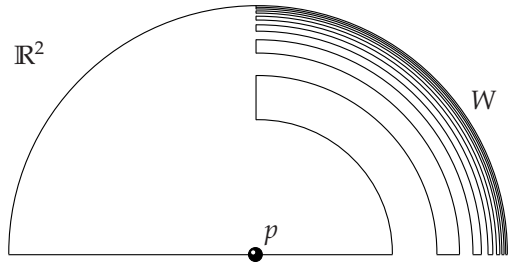

Figure 2: Embedding of the Warsaw circle $W$ into $\mathbb{R}^2$, with the source point $p$. The Reeb graph of the distance function $d(p, \cdot)$ is a circle, whereas $W$ has no simple closed curves.

## 3. Cycle Rank of a Finite Reeb Graph

The *cycle rank*, called also *number of cycles* or *number of loops*, of the Reeb graph $R_f$ is its first Betti number $b_1(R_f)$. We will be interested in those functions whose Reeb graph $R_f$ is a finite topological graph, as is the case for important classes of functions, such as smooth functions with isolated connected components of the critical set (including Morse functions, as well as Morse-Bott functions and functions with isolated critical points) on compact manifolds with possible boundary.

The following bound first appeared in [20, Theorem 13] for any smooth function but only on orientable manifolds:

**Theorem 3.1.** *Let $X$ be a connected locally path-connected topological space and $f : X \to \mathbb{R}$ a continuous function whose Reeb graph $R_f$ is a finite topological graph. Then for the cycle rank $b_1(R_f)$ it holds*

$$b_1(R_f) \le b'_1(X), \tag{8}$$

*where $b'_1(X) = \operatorname{corank}(\pi_1(X))$.*

*Proof.* For a quotient map $\varphi \colon (X, x_0) \to (Y, y_0)$ of topological spaces, where $X$ is locally path-connected and $Y$ is semilocally simply-connected, Calcut et al. [7, Theorem 1.1] showed that if each fiber $\varphi^{-1}(y)$ is connected, then the induced homomorphism $\varphi_* \colon \pi_1(X, x_0) \to \pi_1(Y, y_0)$ is surjective. A connected locally path-connected space is path-connected, so $\pi_1(X)$ is defined unambiguously. For $Y = R_f$, we obtain $\operatorname{corank}(\pi_1(X)) \ge \operatorname{rank}(\pi_1(R_f)) = b_1(R_f)$ since $\pi_1(R_f)$ is a free group. □

Note that in Section 2.2 we defined $b'_1(X)$ only for path-connected spaces. If we relax this requirement in the definition of $b'_1(X)$ by, say, summation by the path-connected components, then the requirement for the space to be connected can be relaxed in Theorem 3.1.

On the other hand, the bound (8) implies a bound

$$b_1(R_f) \le b_1(X), \tag{9}$$

which is widely believed to be obvious for any continuous function on any topological space, as implied, e.g., in [13, §VI.4, p. 141] or [9, Eq. (1)]. However, as Examples 3.3 and 3.5 below show, the requirement for the space to be locally path-connected is essential for the bounds (8) and even (9) to hold:

**Proposition 3.2.** *Let $r, b$ be any non-negative integers. There exist a topological space $X$ and a continuous function $f : X \to \mathbb{R}$ such that for the cycle rank of the Reeb graph $R_f$ and the first Betti number of $X$ it holds*

$$b_1(R_f) = r, \quad b_1(X) = b.$$

The space can be chosen to be a path-connected metric space $(X, d)$ and the function can be chosen as the distance function $f = d(p, \cdot)$ for some $p \in X$ defining the Reeb graph $R_f$ that is a finite topological graph. We will need these properties in Section 5.2. Note that in our examples it holds $b'_1(X) = b_1(X)$.

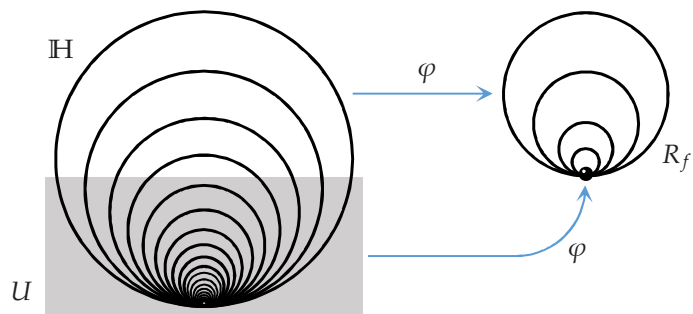

Figure 3: Left: the Hawaiian earring $\mathbb{H}$ with a continuous function $f: \mathbb{H} \to \mathbb{R}$ that is constant in the grayed area $U$ and is the height function outside $U$; $\varphi$ is its Reeb quotient map. Right: the Reeb graph $R_f$ is a finite topological graph.

**Example 3.3.** *Consider a Warsaw circle $W$ as a metric subspace of $\mathbb{R}^2$ shown in Figure 2, with $p$ at the origin and $f = d(p, \cdot)$. It is easy to see that $b_1(W) = 0$ whereas $R_f$ is a circle: $b_1(R_f) = 1$, which contradicts (9). Note that $W$ is not locally connected.*

**Example 3.4.** *Consider a metric subspace $\Theta = \{x_1^2 + x_2^2 = 1\} \cup [p, (1, 0)] \subset \mathbb{R}^2$ with $p$ at the origin; then $b_1(\Theta) = 1$. For the distance function $f = d(p, \cdot)$, $R_f$ is a segment: $b_1(R_f) = 0$, with inequality in both (8) and (9).*

**Example 3.5.** *Given arbitrary $r, b \in \mathbb{N}$, consider the wedge sum $X$ of $r$ copies of $W$ and $b$ copies of $\Theta$ from the examples above, joined by the source point $p$. Then $b_1(R_f) = r$ and $b_1(X) = b$, where $f = d(p, \cdot)$ is the distance function. Note that $X$ is a compact path-connected metric space and the Reeb graph $R_f$ is a finite topological graph.*

The effect might disappear and the requirement for the space to be locally path-connected might be redundant for the bound (9) to hold in a different homology theory, such as Čech homology, related to the shape theory, but we did not check this.

On the other hand, the requirement for the space to be locally path-connected is sufficient but not necessary for (8) to hold for the distance function:

**Example 3.6.** *Consider a path-connected subset $S \subseteq \{x \in \mathbb{R}^n \mid \sum x_i^2 = 1\}$, a point $s \in S$, and $X = S \cup [0, s]$ as a metric subspace of $\mathbb{R}^n$. As in Example 3.4, for $f = d(0, \cdot)$ the Reeb graph $R_f$ is a segment, even though $X$ can be a very ill-behaved space, such as the Warsaw circle from Example 3.3.*

While inequality in (9) holds for many manifolds—see, e.g., Example 2.2,—for a typical continuous function on an orientable surface, such as simple Morse function [9, Lemma A], or even for a wider class of Morse functions [20, Theorem 9], the bound (8) turns into equality. Still inequality in (8) can hold for the distance function on a quite nice space:

**Example 3.7.** *Consider a torus $T^2 = \mathbb{R}^2 / \mathbb{Z}^2$ with the metric induced from $\mathbb{R}^2$. For the distance function $f = d(0, \cdot)$, we have $0 = b_1(R_f) < b_1'(T^2) = 1$.*

Theorem 3.1 can also be used as a bound on $b_1'(X)$:

**Example 3.8.** *The Hawaiian earring $\mathbb{H}$ is a union of countably many decreasing circles. Its fundamental group $\pi_1(\mathbb{H})$ is uncountable and not free [11]. A torsion-free group, such as $\mathbb{Q}$, can have zero co-rank. However, $\mathrm{corank}(\pi_1(\mathbb{H}))$ is at least countably infinite. Indeed, consider a continuous function $f: \mathbb{H} \to \mathbb{R}$ that is constant in a neighborhood $U$ of the "bad" point and is a height function outside $U$; see Figure 3. Then $R_f$ is a finite topological graph; thus (8) gives $b_1'(\mathbb{H}) \geq b_1(R_f)$. The latter can be made arbitrarily large by choosing small enough $U$.*

For connected smooth closed manifolds, the bound (8) is tight even in the class of simple Morse functions:

**Proposition 3.9.** *Let M be a connected smooth closed manifold. There exists a simple Morse function $f : M \to \mathbb{R}$ defining the Reeb graph $R_f$ with the cycle rank*

$$b_1(R_f) = b'_1(M).$$

*Proof.* Since $M$ is compact, the fundamental group $\pi_1(M)$ is finitely generated, so $b'_1(M)$ is finite. By Proposition 2.3, there exist $k = b'_1(M)$ disjoint compact connected codimension-one submanifolds $N_1, \ldots, N_k \subset M$ such that each $N_i$ has a product neighborhood and $M \setminus \bigcup_{i=1}^{k} N_i$ is connected.

By Lemma 2.1, there exists a simple Morse function $f : M \to \mathbb{R}$ such that all $N_i$ are its contours. Since $M \setminus \bigcup_{i=1}^{k} N_i$ is connected, the number of cuts of the Reeb graph $R_f$ is at least $k$, so $b_1(R_f) \geq k$. On the other hand, Theorem 3.1 gives $b_1(R_f) \leq k$. Thus $b_1(R_f) = b'_1(M)$. □

For the case of $\dim M \geq 3$, Michalak [32, Theorem 4.7, Corollary 5.4] has recently shown a stronger fact: in the conditions of Proposition 3.9, for every $b = 0, \ldots, b'_1(M)$ there exists a said function with $b_1(R_f) = b$.

Michalak [33] introduced the notion of the Reeb number $\mathcal{R}(M)$ of a smooth closed manifold $M$ as the maximum cycle rank (he uses the term *number of cycles*) among all Reeb graphs of smooth functions on $M$ with finitely many critical points. As we have just seen (see also [32, Corollary 5.4]), it always holds $\mathcal{R}(M) = b'_1(M)$. However, one can extend the notion of the Reeb number to an arbitrary topological space as the maximum cycle rank $\mathcal{R}(X)$ among all Reeb graphs of continuous functions $f$ on $X$ with the Reeb graph $R_f$ being a finite topological graph. Example 3.5 shows a connected but not locally path-connected space with $\mathcal{R}(X) > b'_1(X)$. Example 3.8 shows a compact connected locally path-connected space $\mathbb{H}$ with $\mathcal{R}(\mathbb{H}) = \infty$.

## 4. Calculation of the Co-rank of the Fundamental Group

Our results involve the value $b'_1(X) = \text{corank}(\pi_1(X))$. Unlike rank, the co-rank is known to be algorithmically computable for finitely presented groups [25, 38]. While we are not aware of any simple method of finding $b'_1(X)$ for an arbitrary space, it can be easily calculated for spaces constructed from simpler ones:[2]

**Theorem 4.1.** *Let $X_1, X_2$ be path-connected topological spaces with finitely generated fundamental groups $\pi_1(X_i)$. Then for the direct product, the union, the wedge sum, and the connected sum it holds:*

 (i) $b'_1(X_1 \times X_2) = \max\{b'_1(X_1), b'_1(X_2)\}$;
 (ii) $b'_1(X_1 \cup X_2) = b'_1(X_1) + b'_1(X_2)$ *if both $X_i$ are open in $X_1 \cup X_2$ and $\emptyset \neq X_1 \cap X_2$ is simply connected;*
 (iii) $b'_1(X_1 \vee X_2) = b'_1(X_1) + b'_1(X_2)$ *if the based topological spaces $(X_i, p_i)$ are locally contractible at $p_i \in X_i$;*
 (iv) $b'_1(M_1 \# M_2) = b'_1(M_1) + b'_1(M_2)$ *for $M_1, M_2$ being closed orientable surfaces[3] or connected n-manifolds with (possibly empty) boundary, $n \geq 3$, with finitely generated fundamental group (e.g., compact manifolds).*

*Proof.* (i) follows from $\pi_1(X_1 \times X_2) \cong \pi_1(X_1) \times \pi_1(X_2)$ along the lines of [18, Theorem 3.1]. (ii) The Seifert–van Kampen theorem gives $\pi_1(X_1 \cup X_2) \cong \pi_1(X_1) * \pi_1(X_2)$, while for finitely generated groups, [24, Proposition 6.4] states $\text{corank}(G_1 * G_2) = \text{corank}(G_1) + \text{corank}(G_2)$. (iii) can be reduced to (ii) by extending each $X_i \subset X_1 \vee X_2$ by a small neighborhood of $p$, where $p$ denotes the identified points $p_1$ and $p_2$. (iv) is similar to (ii) for $n \geq 3$; for closed orientable surfaces $M_g^2$ of genus $g$, it is verified directly since $b'_1(M_g^2) = g$ [23]. □

For example, the values of $b'_1$ given in Example 2.2 for $T^n = \times^n S^1$ and $M_g^2 = \#^g T^2$ can be calculated using this theorem. More details on the calculation of $b'_1(M)$ for different manifolds can be found in [18].

In Theorem 4.1 (i) (and Theorem 4.3 (i) below), the direct product cannot be replaced with arbitrary fiber bundle: for instance, $S^3$ is a fiber bundle of $S^1$ over $S^2$. Example 2.2 shows that the conditions in Theorem 4.1 (iv) (and Theorem 4.3 (iv) below) for the surface to be closed and orientable are essential.

---

[2] Some of the facts given in this section have been stated for smooth manifolds, and some only for the orientable case, in [18, Theorem 3.1], [19, Theorems 21, 27], and [24, Proposition 6.4]. Here we extend them to more general topological spaces, and also consider a union and wedge sum of spaces. The facts about manifolds are given here for completeness.

[3] For the case of non-orientable surfaces and surfaces with boundary, see Example 2.2.

While there exist path-connected topological spaces with finitely-generated fundamental group that are not locally contractible at a point (e.g., a cone over the Hawaiian earring), we currently do not know whether the condition for the spaces to be locally contractible at the points $p_i$ in Theorem 4.1 (iii) is essential.

In the real world (and in computational geometry), shapes often have variable dimension: a table with legs or a cup with a handle. Also, there are joining operations more general than wedge or connected sum: a cup is joined to the table by its bottom. Theorem 4.1 can be generalized to a wider range of operations:

**Remark 4.2.** *One can generalize the wedge sum to gluing by a larger subset with corresponding conditions on (local) contractibility. One also can generalize the connected sum to arbitrary topological spaces with a neighborhood homeomorphic to $\mathbb{R}^n$. Then Theorem 4.1 will hold for such operations.*

For some spaces $X$, the value $b'_1(X)$ is bounded from above by the isotropy index $h(X)$, which is easier to calculate; see Section 2.3. In turn, we can derive for $h(X)$ equations similar to those from Theorem 4.1:

**Theorem 4.3.** *Under the same conditions as in Theorem 4.1 on each equation, it holds*

(i)  $h(X_1 \times X_2) = \max\{h(X_1), h(X_2)\}$,
(ii) $h(X_1 \cup X_2) = h(X_1) + h(X_2)$,
(iii) $h(X_1 \vee X_2) = h(X_1) + h(X_2)$,
(iv) $h(M_1 \# M_2) = h(M_1) + h(M_2)$.

*For (ii), instead of $X_1 \cap X_2$ being simply connected, it is enough to require $H^1(X_1 \cap X_2; \mathbb{Q}) = 0$.*

*Proof.* By (4), (i) can be obtained along the lines of [19, Theorem 27]. For (ii), since $X_1 \cap X_2$ is simply connected, we have $H^1(X_1 \cap X_2; \mathbb{Q}) = 0$. The Mayer–Vietoris sequence gives $H^1(X_1 \cup X_2; \mathbb{Q}) = H^1(X_1; \mathbb{Q}) \oplus H^1(X_2; \mathbb{Q})$; the cup product translates into component-wise product. Applying [19, Lemma 20] to these vector spaces gives (ii). Now (iii) and (iv), as in Theorem 4.1, can be reduced to (ii). For the orientable surfaces $M_g^2$, (iv) is verified directly since $h(M_g^2) = g$; see, e.g., Proposition 2.7. □

Given the similarity between Theorems 4.1 and 4.3, for many classical manifolds, as well as their connected sums and the directed products, it holds $b'_1(M) = h(M)$, though strict inequality in (3) is also possible [19, Example 41]. Theorem 4.3 can also be generalized along the lines of Remark 4.2.

## 5. Approximation of Spaces by Graphs and Bounds on the Distortion of the Reeb Quotient Map

Theorem 3.1 can be used to improve bounds on the approximation of some spaces by simple enough graphs, and in particular by Reeb graphs, as well as on the distortion of the Reeb quotient map. We will apply it to the distance function on a geodesic space and to a simple Morse function on a Riemannian manifold (we assume all Riemannian manifolds to be smooth), where our bound (8) provides considerably better estimates than the trivial bound (9) used by previous authors.

### 5.1. Approximation of a Compact Geodesic Space by Metric Graphs

Let $(X, d)$ be a compact geodesic space such that $b_1(X)$ is finite, and the source point $p \in X$ be such that the Reeb graph $R$ defined by the distance function $d(p, \cdot)$ is finite. It is known that

$$\delta_i^X = \inf\{d_{GH}(X, G) \mid G \text{ is a finite metric graph with } b_1(G) \leq i\} \to 0 \tag{10}$$

when $i \to \infty$ [6, Proposition 7.5.5]. However, for practical applications, such as those in computational geometry [4], it is very important to know how fast this sequence starts to converge: whether $X$ can be approximated well by graphs simple enough for real-life computation. Studying the rate of convergence of this sequence with small $i$, Mémoli and Okutan [31, Theorem 1.2] found the bounds

$$\frac{1}{16i + 12}\rho \leq \delta_i^X \leq \rho, \qquad \text{if } i \geq b, \tag{11}$$

$$\frac{1}{16b + 12}\rho \leq \delta_i^X \leq \rho + 6(b+1)a_{i+1}^X \qquad \text{if } i < b, \tag{12}$$

where $\rho = d_{GH}(X, R)$, the $a_1^X \geq a_2^X \geq \ldots$ are the lengths of the intervals in the first persistent barcode of the open Vietoris-Rips filtration of $X$ [31, Section A.2], and effectively $b = b_1(R)$. However, in terms of the characteristics of the space they report their results with $b = b_1(X)$ using the trivial bound (9). Changing their $b = b_1(X)$ to $b = b'_1(X)$ implied by our bound (8) improves this estimate:

**Example 5.1.** *For n-torus $T^n$, with $b = b_1(T^n) = n$ and $b'_1(T^n) = 1$, for any $i > 0$ our bound gives (11), while [31] uses for $i < n$ the much weaker estimate (12): for instance, for $i = 1$, we obtain much tighter estimate $\frac{1}{28}\rho \leq \delta_i^X \leq \rho$ than their $\frac{1}{16n+12}\rho \leq \delta_i^X \leq \rho + 6(n+1)a_2^X$. This reflects the intuition, missed by the expression as reported in [31, Theorem 1.2], that even a very high-dimensional torus can be approximated by a circle-like graph.*

### 5.2. Bound on the Distortion: The Case of the Distance Function

Chazal et al. [8, Propositions 1 and 2] considered a compact geodesic space $(X, d)$ and a function $f : X \to \mathbb{R}$ that is the distance $f = d(p, \cdot)$ from a source point $p \in X$. For the distortion of the Reeb quotient map $\varphi$, they showed that if the Reeb graph $R_f$ is a finite topological graph, then

$$\text{dis}(\varphi) \leq 2(b_1(R_f) + 1)D, \tag{13}$$

where $D$ is the supremum of the diameters of contours of $f$.

As to estimating $\text{dis}(\varphi)$ in terms of topological characteristics of the space $X$ itself, they only mentioned the bound (9), which does hold for geodesic spaces since they are locally path-connected [1, Corollary 2.1]. However, our Theorem 3.1 applied to (13) gives a much stronger bound:

**Theorem 5.2.** *Let $(X, d)$ be a compact geodesic space, $p \in X$, and $f : M \to \mathbb{R}$ the distance function $f = d(p, \cdot)$. Let the Reeb graph $R_f$ be a finite topological graph. Then for the distortion of the Reeb quotient map $\varphi$ it holds*

$$\text{dis}(\varphi) \leq 2(b'_1(X) + 1)D,$$

*where $b'_1(X) = \text{corank}(\pi_1(X))$ and $D$ is the supremum of the diameters of contours of $f$. In particular, for the Gromov-Hausdorff distance, it holds*

$$d_{GH}(X, R_f) \leq (b'_1(X) + 1)D.$$

**Example 5.3.** *For n-torus $T^n$, our bound gives $\text{dis}(\varphi) \leq 4D$, see Example 2.2, whereas the bound (9) referred to by Chazal et al. [8] only gives $\text{dis}(\varphi) \leq 2(n+1)D$. For the closed surface $M_g^2$ of genus $g$, we obtain $\text{dis}(\varphi) \leq 2(g+1)D$, whereas (9) gives $\text{dis}(\varphi) \leq 2(2g+1)D$.*

It is essential both for the bound (9) and for our stronger bound (8) that the considered space is locally path-connected. As Example 3.5 shows, otherwise the homology of the space $X$ does not provide any useful information on the cycle rank of $R_f$ even for the distance function.

### 5.3. Bound on the Distortion: The Case of a Simple Morse Function

Let $M$ be a closed connected $n$-dimensional Riemannian manifold, $n \geq 2$, with the distance function $d(\cdot, \cdot)$, and $f : M \to \mathbb{R}$ an $L$-Lipschitz simple Morse function. Consider $p \in M$ and denote

$$B(b) = 4(b+1)^2 \left( \left( \frac{2L}{b+1} \cdot \frac{\mu^n(M)}{T_f} \right)^{\frac{1}{n}} + 8 \left( (\text{diam}(M))^{\frac{1}{n}} \varepsilon_p^{\frac{n-1}{n}} + \varepsilon_p \right) \right) + |L - 1| \, \text{diam}(M),$$

where $T_f$ is the thickness of $f$ described in Section 2.5, $\mu^n$ is the $n$-dimensional Hausdorff measure on $M$, and $\varepsilon_p = \|f - d(p, \cdot)\|_\infty$.

**Theorem 5.4.** *Let $R_f$ be the Reeb graph of $f$ and $\varphi: M \to R_f$ the Reeb quotient map. Then for the metric distortion of $\varphi$ it holds:*

$$\mathrm{dis}(\varphi) \leq B(b'_1(M)), \tag{14}$$

*where $b'_1(M) = \mathrm{corank}(\pi_1(M))$. In particular, for the Gromov-Hausdorff distance, it holds*

$$d_{GH}(X, R_f) \leq \frac{1}{2} B(b'_1(M)).$$

Note that under the conditions of the theorem the Reeb graph $R_f$ is a finite topological graph [39, Theorem 1].

*Proof.* The proof of [30, Proposition 4.4] implies that, similarly to (13),

$$\mathrm{dis}(\varphi) \leq (2b_1(R_f) + 1)(D + 4\varepsilon_p) + |L - 1| \mathrm{diam}(M), \tag{15}$$

where $D$ is the supremum of the diameters of contours of the function $f$. In turn, the proof of [30, Proposition 6.2] implies a bound

$$D \leq \left(2^{n+1} k^{n-1} L \cdot \frac{\mu^n(M)}{T_f}\right)^{\frac{1}{n}} + 8k \left((\mathrm{diam}(M))^{\frac{1}{n}} \varepsilon_p^{\frac{n-1}{n}} + \varepsilon_p\right), \tag{16}$$

where $k$ is the minimum number such that for any $k$ regular contours $N_1, \ldots, N_k$ of $f$, the set $M \setminus \bigcup_{i=1}^k N_i$ is disconnected. Applying Theorem 3.1 to $b_1(R_f)$ in (15) and Remark 2.4 to $k$ in (16), we obtain (14). □

Mémoli and Okutan [30, Theorem 1.1] used coarser bounds $b_1(R_f) \leq b_1(M)$ and $k \leq b_1(M) + 1$ [30, Lemma 2.7] and obtained a bound

$$\mathrm{dis}(\varphi) \leq B(b_1(M)). \tag{17}$$

Since $b'_1(M) \leq b_1(M)$ and $B(b)$ grows with $b$, our bound (14) implies (17), in many cases being much tighter. Indeed, simple Morse functions are dense in the space of continuous functions on a closed manifold. Consider distance-like simple Morse functions, i.e., assume small $\varepsilon_p$ and thus $L$ close to 1. Denoting Mémoli and Okutan's bound by $B = B(b_1(M))$ and our bound by $B' = B(b'_1(M))$, both with $\varepsilon_p = 0$ and $L = 1$, for such functions we have

$$\frac{B}{B'} = \left(\frac{b_1(M) + 1}{b'_1(M) + 1}\right)^{2 - \frac{1}{n}}.$$

Example 2.2 and Section 4 show that this ratio is often quite significant:

**Example 5.5.** *For an n-torus, $b_1(T^n) = n$, while $b'_1(T^n) = 1$, which gives $\frac{B}{B'} = \left(\frac{n+1}{2}\right)^{2 - \frac{1}{n}}$. For $n \geq 3$, our bound is more than twice tighter; for $n \gg 1$, the improvement is by almost $\frac{1}{4} n^2$.*

**Example 5.6.** *For the closed surface $M_g^2$ of genus $g$, we have $b_1(M_g^2) = 2g$, while $b'_1(M_g^2) = g$, so $\frac{B}{B'} = (2 - \frac{1}{g+1})^{\frac{3}{2}}$. For $g \geq 2$, our bound is more than twice tighter; for $g \gg 1$, the improvement is by almost $2^{\frac{3}{2}} \approx 2.8$.*

For $M_g^2$, substituting 6 for $4\sqrt{2}$ to compensate for $\varepsilon_p$ and $L$, our bound

$$\mathrm{dis}(\varphi) \leq 6(g+1)^{\frac{3}{2}} \sqrt{\mu^2(M_g^2)}$$

coincides with the one previously obtained by Zinov'iev [43, Proposition 3.1].

*5.4. Bounds from Below*

Upon discussing the upper bounds on the distortion $\mathrm{dis}(\varphi)$ for the Reeb quotient map, one can wonder how well a metric space can be approximated by the Reeb graphs of continuous functions.

Some spaces can be approximated arbitrarily well even by Reeb graphs that are finite graphs:

**Example 5.7.** *In the following examples, we assume that the spaces are equipped with the intrinsic metric obtained from the metric induced by a suitable embedding of the space into $\mathbb{R}^n$.*

*For a finite topological graph $G$, $p \in G$, and $f(\cdot) = d(p, \cdot)$, it holds $R_f \cong G$, i.e., $\mathrm{dis}(\varphi) = 0$.*

*For the Hawaiian earring $\mathbb{H}$ and the function $f$ being the distance function from a ball $B_r$ (grayed area in Figure 3) of radius $r$ around the "bad" point $p$, we have $\mathrm{dis}(\varphi) \neq 0$ since $R_f \neq \mathbb{H}$, but $\mathrm{dis}(\varphi) \to 0$ when $r \to 0$. Note that for $f(\cdot) = d(p, \cdot)$, we have $\mathrm{dis}(\varphi) = 0$, but $R_f$ is not a finite graph.*

By (10), all good enough spaces can be arbitrarily well approximated by finite graphs; however, not in all cases by Reeb graphs. Considering only continuous functions whose Reeb graphs are finite topological graphs (which are one-dimensional simplicial complexes), the distortion of the Reeb quotient map $\varphi \colon X \to R_f$ is bounded by

$$\mathrm{dis}(\varphi) \geq u_1(X),$$

the one-dimensional Urysohn width (5), which in turn in some cases can be bounded by characteristics of the space, e.g., the volume in the case of a Riemannian manifold; see (7). In the general case of a continuous function $f$, it holds

$$\mathrm{dis}(\varphi) \geq \mathrm{diam}(C)$$

for any contour $C$ of $f$, so it is enough to bound the diameter of the largest contour of $f$.

**Proposition 5.8.** *Let $X$ be a metric space containing a homeomorphic (not necessarily isometric) image of $\mathbb{R}^2$ (e.g., an $n$-dimensional Riemannian manifold, $n \geq 2$). Then there exists a value $b(X) > 0$ such that any continuous function $f \colon X \to \mathbb{R}$ has a contour $C$ with*

$$\mathrm{diam}(C) \geq b(X).$$

*Proof.* Consider $U \subseteq X$ with $\psi \colon \mathbb{R}^2 \to U$ being a homeomorphism, and a closed Euclidean ball $E \subseteq \mathbb{R}^2$, which we can assume to be a unit ball. Denote $S = \{(x, y) \in E^n \times E^n \mid |x - y| \geq \sqrt{3}\}$. Since $S$ is compact, there exists $b = \min_S d(\psi(x), \psi(y))$, where $d(\cdot, \cdot)$ is the distance in $X$, with $b > 0$ since $x \neq y$ for all $(x, y) \in S$.

Let $f \colon X \to \mathbb{R}$ be continuous. By Proposition 2.8, $F = f \circ \psi$ has a contour $L$ containing two points $(x, y) \in S$. Then $\psi(L)$, a subset of a contour $C$ of $f$, contains points $p = \psi(x)$, $q = \psi(y)$ with $d(p, q) \geq b$. □

In Section 6, we give more specific bounds for Riemannian manifolds and discuss related concepts.

## 6. Diameters of Contours on Disks and on Riemannian Manifolds

In Section 5.4, we have shown that some metric spaces $X$ cannot be arbitrarily well approximated by Reeb graphs: for the distortion of the Reeb quotient map, it holds

$$\inf_{f \colon X \to \mathbb{R}} \mathrm{dis}(\varphi) \geq b(X),$$

with $b(X) > 0$, for example, for spaces containing a homeomorphic image of $\mathbb{R}^2$. For the latter value, we can introduce the following notation:

**Definition 6.1.** *The* Reeb width *of a metric space $X$ is the infimum, over all real-valued continuous functions on $X$, of the suprema of the diameters of their contours:*

$$b(X) = \inf_{\substack{\text{continuous} \\ f \colon X \to \mathbb{R}}} \sup_{\text{contours } C \text{ of } f} \mathrm{diam}(C).$$

For example, for Proposition 2.8 implies $b(\mathbb{R}^n) = \infty$, $n \geq 2$.

Our definition of the Reeb width $b(X)$ closely resembles that of Urysohn width (5). Indeed, when appropriate, one can consider the Reeb graph $R_f$ as $Y^1$ and the Reeb quotient map as $\varphi$ in (5). According to [26, Corollary 3.9], for a unit disk $D \subset R^2$, it holds $b(D) = u_1(D) = \sqrt{3}$. For the Hawaiian earring, we also have $b(\mathbb{H}) = u_1(\mathbb{H}) = 0$: while for each $Y$, the diameter of some level sets is positive, they can be made arbitrarily small; see Figure 3. We leave it for the future work to determine whether and when the two notions are different.

For a Riemannian manifold, the Reeb width can be bounded from below in terms of its characteristics; see Section 2.6 for notation and details.

First, we will need a refined version of Proposition 2.8, which due to its generality has value in itself:

**Theorem 6.2.** *Let $E^n \subset \mathbb{R}^n$, $n \geq 2$, be a closed unit ball and $f : E^n \to \mathbb{R}$ a continuous function. Then $f$ has a contour $C$ such that*

  (i) *$n = 2$: $\mathrm{diam}(C \cap \partial E) \geq \sqrt{3}$,*
  (ii) *$n \geq 3$: $\mathrm{diam}(C \cap \partial E) = 2$.*

*Proof.* The fact about $\mathrm{diam}(C)$ is stated in Proposition 2.8, so we only need to show that two points $x, y \in C$ with the corresponding distance $\|x - y\|$ exist *at the boundary* of the ball. The case $n \geq 3$ is trivial, so we assume $E = E^2$ to be a unit disk around the origin 0.

Denote by $R(t_1, t_2)$ the ring-shaped area $\{t_1 \leq \|p\| \leq t_2\} \subset E$ for given $0 < t_1 < t_2 \leq 1$. For each $k = 1, 2, \ldots$, consider a unit ball $E'$ around $0'$ and a homeomorphism $\psi_k : E' \to E$ that stretches $E'$ near $0'$ and shrinks it near $\partial E'$ along radii (preserving angles at the origin), mapping the respective ring-shaped area $R'(\sqrt{3} - 1, 1) \subset E'$ to $R(1 - \frac{1}{k}, 1) \subset E$. Consider the function $f'_k = f \circ \psi_k : E' \to \mathbb{R}$. By Proposition 2.8, $f'_k$ has a contour $C'_k$ with two points $x'_k, y'_k \in C'$ such that $\|x'_k - y'_k\| \geq \sqrt{3}$; thus $x'_k, y'_k \in R'_{2-\sqrt{3}}$ and $\angle x'_k 0' y'_k \geq \frac{2\pi}{3}$. Then $C_k = \psi_k(C'_k)$ is a contour of $f$ and for $x_k = \psi_k(x'_k)$, $y_k = \psi_k(y'_k)$ we have $x_k, y_k \in C \cap R(1 - \frac{1}{k}, 1) \subset E$, with $\angle x_k 0 y_k \geq \frac{2\pi}{3}$ and thus $\|x_k - y_k\| \geq \sqrt{3}$.

We have built a sequence of pairs $(x_k, y_k)$, each lying on a contour $C_k$ of $f$, with $\|x_k - y_k\| \geq \sqrt{3}$, located at the distance of $\frac{1}{k}$ from the boundary $\partial E$. Since $E$ is compact, there are points $x, y \in E$ such that $x_k \to x$, $y_k \to y$ (up to selection of a subsequence). Obviously, $x, y \in \partial E$, $\|x - y\| \geq \sqrt{3}$, and $f(x) = f(y)$; without loss of generality, assume $f(x) = f(y) = 0$.

Now we only need to show that $x, y$ lie in the same contour of $f$. For each $i = 1, 2, \ldots$, consider $S_i = f^{-1}[-\frac{1}{i}, \frac{1}{i}]$, with $x, y \in \mathrm{Int}\, S_i$. Consider small balls $X_i, Y_i \subseteq S_i$ around $x$ and $y$, respectively. Up to a finite number, we have $x_k \in X_i$ and $y_k \in Y_i$; in particular, $x$ and $y$ belong to the same connected component of $S_i$, which we denote by $K_i$; it is closed since $S_i$ is closed. Note that $K_{i+1} \subseteq K_i$.

We obtained a nested sequence of closed non-empty connected subsets of a compact metric space $E$; it is known that $K = \bigcap_{i=1}^{\infty} K_i$ is connected. Since $f(K) = 0$ and $x, y \in K$, they have the desired properties. □

While for $n \geq 3$ there actually exists a contour $C$ and points $x, y \in C \cap \partial E$ with $\|x - y\| = 2$, in the case of $n = 2$ there may be no contour $C$ with $C \cap \partial E$ containing points with $\|x - y\| = \sqrt{3}$, as the example of the function $f(x^1, x^2) = |x^1| + |x^2|$ shows.

As a side note, using $R(1 - \frac{1}{k}, 1 - \frac{1}{2k}) \subset E$ instead of $R(1 - \frac{1}{k}, 1) \subset E$ in the proof above gives a variant of the statement for an open ball:

**Proposition 6.3.** *For an open (or closed) $n$-dimensional unit ball $E \subset \mathbb{R}^n$, $n \geq 2$, it holds:*

  (i) *$n = 2$: $b(E) \geq \sqrt{3}$,*
  (ii) *$n \geq 3$: $b(E) = 2$,*

*and for any continuous function $f : E \to \mathbb{R}$ there exist $x, y \in \partial E$ with $\|x - y\| \geq \sqrt{3}$ ($n = 2$) or $\|x - y\| = 2$ ($n \geq 3$) such that for any neighborhoods $U = U(x)$, $V = V(y)$ there exists a contour $C$ of $f$ with $C \cap U \neq \emptyset$ and $C \cap V \neq \emptyset$.*

Note that the function $f(x) = \|x\|$ on an open ball, $n \geq 3$, does not have any contour of diameter 2. We can conjecture that for an open two-dimensional Euclidean unit ball $\mathring{E}$, any continuous function $f : \mathring{E} \to \mathbb{R}$ has a contour $C$ with $\operatorname{diam}(C) \geq \sqrt{3}$, but we do not have a proof of it.

**Proposition 6.4.** *Let $M$ be an $n$-dimensional Riemannian manifold, $n \geq 2$. Let $p \in M$, $r \in \mathbb{R}$ be such that $0 < r \leq \min\{\operatorname{conv.rad}(p), \operatorname{inj}(p)\}$,[4] $B \subset M$ be a closed ball of radius $r$ around $p$, and $K = \sup_{x \in B} \sec(x)$, where $\sec(x)$ is the sectional curvature. Then for the Reeb width of $M$, it holds:*

$$b(M) \geq \begin{cases} 2r, & n \geq 3, \quad (18) \\ \sqrt{3}r, & n = 2, K \leq 0, \quad (19) \\ \frac{2}{\sqrt{K}} \arcsin\left(\frac{\sqrt{3}}{2} \sin\left(r\sqrt{K}\right)\right), & n = 2, K > 0, r < \frac{\pi}{2\sqrt{K}}, \quad (20) \\ \frac{2\pi}{3\sqrt{K}}, & n = 2, K > 0, r \geq \frac{\pi}{2\sqrt{K}}. \quad (21) \end{cases}$$

The cases 20 and 21 represent $\frac{2}{\sqrt{K}} \arcsin\left(\frac{\sqrt{3}}{2} \sin \min\left\{r\sqrt{K}, \frac{\pi}{2}\right\}\right)$.

*Proof.* Let $f : M \to \mathbb{R}$ be a continuous function. As in the proof of Proposition 5.8, consider a closed Euclidean ball $E \subset T_p M$ of radius $r$ around the origin. Since $r$ is finite, $\partial E \neq \emptyset$; since $r \leq \operatorname{inj}(p)$, the exponential map $\exp_p : E \to B$ is a diffeomorphism, and since $r \leq \operatorname{conv.rad}(p)$, distances $d(\cdot, \cdot)$ in $B$ are measured along geodesic segments lying in it. Consider the function $F = f \circ \exp_p : E \to \mathbb{R}$.

(18): By Proposition 2.8, some contour $L \subseteq E$ of $F$ contains points $x, y \in E$ with $\|x - y\| = 2r$, thus they are antipodal points, i.e., the straight segment $[x, y]$ passes through the origin. Then $\exp_p(L)$, a part of a contour $C$ of $f$, contains the points $\exp_p(x), \exp_p(y)$ connected by the geodesic segment $\exp_p([x, y])$ of length $2r$, which is the shortest path between them by the choice of $B$. We obtain $\operatorname{diam} C \geq 2r$.

(19): Similarly, by Proposition 2.8, a contour $L$ of $F$ contains points $x, y$ with $\|x - y\| \geq \sqrt{3}r$. Taking $T_pM$ as the model space $S_0$ in Theorem 2.9, we conclude as above that $\operatorname{diam}(C) \geq d(\exp_p(x), \exp_p(y)) \geq \|x - y\| \geq \sqrt{3}$.

(20): Applying again Theorem 2.9, we reduce the task to calculating the distances in the sphere $S_K$ of curvature $K$ as a model space, which we represent as a sphere in $\mathbb{R}^3$ of radius $\frac{1}{\sqrt{K}}$. Consider $x, y \in \partial E$ given by Theorem 6.2 (i). For two points $q = \exp_p(x)$, $s = \exp_p(y)$ separated by the azimuth angle $\alpha$ at the line of latitude with the polar angle determined by $r$, direct calculation in spherical coordinates gives

$$d(q, s) = \frac{2}{\sqrt{K}} \arcsin\left(\sin \frac{\alpha}{2} \sin\left(r\sqrt{K}\right)\right).$$

The condition $\|x - y\| \geq \sqrt{3}r$ gives $\alpha \geq \frac{2\pi}{3}$, thus $\sin \frac{\alpha}{2} \geq \frac{\sqrt{3}}{2}$.

(21): Since this expression reaches a maximum by $r$ at $r_{\max} = \frac{\pi}{2\sqrt{K}}$, in the case of $r > r_{\max}$ we consider a smaller ball of radius $r_{\max}$ inside $B$. □

As an approximation, the expression (20) can be simplified by replacing $\sin(x)$ with $\frac{2}{\pi}x$ and $\arcsin(x)$ with $x$, which results in $b(M) \geq \frac{2\sqrt{3}}{\pi}r \approx 1.10\,r$, while (21) reduces to $b(M) \geq \frac{2\pi}{3\sqrt{K}} \approx 2.09\frac{1}{\sqrt{K}}$:

**Corollary 6.5.** *Under the conditions of Proposition 6.4, it holds*

$$b(M) \geq \min\left\{r, \frac{2}{\sqrt{K}}\right\}. \quad (22)$$

While the simplified bound 22 is not tight, the bounds (20) and 21 are tight:

---

[4] Generally, both $\operatorname{inj}(p) \ll \operatorname{conv.rad}(p)$ (Example 6.10) and $\operatorname{conv.rad}(p) \ll \operatorname{inj}(p)$ [12, Theorem 3.4] are possible.

**Example 6.6.** *Consider the upper hemisphere of the unit sphere in $\mathbb{R}^3$ as a Riemannian manifold M. Taking its pole as p and $r = \frac{\pi}{2}$ in Proposition 6.4, the bound in (21) gives $b(M) \geq \frac{2\pi}{3}$, since $K = 1$.*

*Denote by Y the tripod-shaped union of three meridians at equal angles from each other. Let $f : M \to \mathbb{R}$ be the distance from a point to Y. Then Y is its contour, with $\operatorname{diam}(Y) = \frac{2\pi}{3}$. Any other contour will lie in one of the connected components of $M \setminus Y$, which are triangles with the same diameter $\frac{2\pi}{3}$. We obtain $b(M) = \frac{2\pi}{3}$.*

Note that in Example 6.6, the bounds (18) and (19) do not apply: with $r = \frac{\pi}{2}$, the bound (19) gives $b(M) \geq \frac{\sqrt{3}\pi}{2} \approx 2.72$, while in fact $b(M) = \frac{2\pi}{3} \approx 2.09$. Also note that Example 6.6 cannot be extended to the whole unit sphere, since its Reeb width is $\pi$ [26, Theorem 2.2].

For a compact Riemannian manifold $M$, Proposition 2.10 allows expressing Proposition 6.4 in terms of global characteristics of the manifold by assuming $r = \min\left\{\frac{1}{2}\operatorname{inj}(M), \frac{\pi}{2\sqrt{K}}\right\}$ and $K = \sec(M)$; recall that $\sec(M) = \sup_{p \in M} \sec(p)$:

**Theorem 6.7.** *Let M a compact Riemannian manifold, $K = \sec(M)$, and*

$$b(M) = \inf_{\substack{\text{continuous}\\ f : M \to \mathbb{R}}} \sup_{\text{contours } C \text{ of } f} \operatorname{diam}(C)$$

*denote the Reeb width of M. Then*

$$b(M) \geq \begin{cases} \operatorname{inj}(M), & \dim M \geq 3, K \leq 0, \\ \min\left\{\operatorname{inj}(M), \frac{\pi}{\sqrt{K}}\right\}, & \dim M \geq 3, K > 0, \\ \frac{\sqrt{3}}{2}\operatorname{inj}(M), & \dim M = 2, K \leq 0, \\ \frac{2}{\sqrt{K}}\arcsin\left(\frac{\sqrt{3}}{2}\sin\left(\frac{\sqrt{K}}{2}\operatorname{inj}(M)\right)\right), & \dim M = 2, K > 0, \operatorname{inj}(M) < \frac{\pi}{\sqrt{K}}. \qquad (23)\\ \frac{2\pi}{3\sqrt{K}}, & \dim M = 2, K > 0, \operatorname{inj}(M) \geq \frac{\pi}{\sqrt{K}}. \qquad (24) \end{cases}$$

The cases (23) and (24) represent $\frac{2}{\sqrt{K}}\arcsin\left(\frac{\sqrt{3}}{2}\sin\min\left\{\frac{\sqrt{K}}{2}\operatorname{inj}(M), \frac{\pi}{2}\right\}\right)$.

Similarly, Corollary 6.5 can be rewritten as the following highly simplified bound; the term $\frac{\pi}{2\sqrt{K}}$ present in Proposition 2.10 (i) is omitted here by (24):

**Corollary 6.8.** *Let M a compact Riemannian manifold and $K = \sec(M)$. Then*

$$b(M) \geq \min\left\{\frac{1}{2}\operatorname{inj}(M), \frac{2}{\sqrt{K}}\right\}.$$

Such bounds, while convenient since they use known characteristics of the manifold, may be misleading as to whether the Reeb width, and thus the distortion $\operatorname{dis}(\varphi)$ of the Reeb quotient map, is determined by the regions with the highest curvature, while in fact it is determined by large regions of small curvature:

**Example 6.9.** *Consider a flattened sphere, or a thick disk, in $\mathbb{R}^3$ (the outer surface of a saucer). Since its curvature at the edge is high, Theorem 6.7 gives a small value for $b(M)$ and thus for $\operatorname{dis}(\varphi)$, while Proposition 6.4 gives high value due to the large flat region in the center of the saucer.*

**Example 6.10.** *For a non-complete manifold, such as an open unit disk in $\mathbb{R}^n$, or for a manifold with (non-empty) boundary, such as a closed unit disk, we have $\operatorname{inj}(M) = 0$, so the bounds in the global terms from Theorem 6.7 are not useful (even though the curvature is bounded, with K=0), while the bounds in the local terms from Proposition 6.4 are exact for the unit disk (or for the manifold from Example 6.6) considering the center as the point p.*

There are other interesting bounds and on conv. rad($M$) and inj($M$); see, e.g., [12, 42], with may result in other versions of Theorem 6.7.

Local versions of Proposition 2.10 have been discussed by the community (though we did not find them in trusted citable sources), with which the requirement of compactness in Theorem 6.7 might be relaxed to, for example, boundedness of the curvature and completeness of the manifold. We leave such generalizations for future work.

In the case of negative curvature, the bounds (18) and (19) can be improved using a triangle comparison theorem and direct calculations on the hyperbolic plane along the lines of (20). However, further refining the bounds on $b(M)$ is outside the scope of this paper and is left for the future work. We also leave as an open question whether there is a non-trivial bound on $d_{GH}(R_f, M)$ in terms of the characteristics of the manifold. We can conjecture that it exists and is closely related to the Reeb width $b(M)$.